%% file: manuscript.tex
\documentclass[reqno]{amsart}

\usepackage[margin=3cm]{geometry}

\usepackage{mystyle}
\usepackage{fontawesome}

\makeatletter
\def\blfootnote{\gdef\@thefnmark{}\@footnotetext}
\makeatother

\title[A general multi-stratum model for a nanofunctionalized
    releasing capsule]{
    A general multi-stratum model for a nanofunctionalized
    releasing capsule: a computational study
}

\author[E.\ Onofri, E.\ Cristiani, A.\ Martelli, P.\ Gentile, J.\ G.\ Hern\'andez, and G.\ Pontrelli]{}

\begin{document}

\blfootnote{$^{\star}$ Corresponding Author, (\href{mailto:elia.onofri@cnr.it}{\faEnvelopeO}) \texttt{elia.onofri@cnr.it}, (\href{https://www.eliaonofri.it}{\faGlobe}) \texttt{www.eliaonofri.it}}

\maketitle

\vspace{-1em}

\begin{center}
    \begin{minipage}{.89\linewidth}\centering
        \textsc{Elia Onofri}$^{a, b, \star, \orcidlink{}}$,
        \textsc{Emiliano Cristiani}$^{a, \orcidlink{0000-0002-7015-2371}}$,
        \textsc{Andrea Martelli}$^{c, \orcidlink{0000-0003-0538-4480}}$,\\
        \textsc{Piergiorgio Gentile}$^{d, e, \orcidlink{0000-0002-3036-6594}}$,
        \textsc{Joel Gir\'on Hern\'andez}$^{f, d, \orcidlink{0000-0003-1245-4475}}$,
        \textsc{Giuseppe Pontrelli}$^{a, \orcidlink{0000-0003-0659-210X}}$.
        \\
        \bigskip
        \begin{minipage}{.45\linewidth}\centering
            \footnotesize
            $^a$Istituto per le Applicazioni del Calcolo (IAC),\\
            Consiglio Nazionale delle Ricerche (CNR)\\
            Rome 00185, Italy
        \end{minipage}
        \hfill
        \begin{minipage}{.45\linewidth}\centering
            \footnotesize
            $^b$Istituto per le Applicazioni del Calcolo (IAC),\\
            Consiglio Nazionale delle Ricerche (CNR)\\
            Naples 80131, Italy
        \end{minipage}
        
        \medskip

        \begin{minipage}{.45\linewidth}\centering
            \footnotesize
            $^c$Dipartmento di Ingegneria, E.\ Ferrari,\\
            Università di Modena e Reggio Emilia\\
            Modena 41121, Italy
        \end{minipage}
        \hfill
        \begin{minipage}{.45\linewidth}\centering
            \footnotesize
            $^d$Centro de Biomateriales e Ingeniería Tisular\\
            Universitat Politècnica de València (CBIT)\\
            València 46022, Spain
        \end{minipage}
        
        \medskip

        \begin{minipage}{.45\linewidth}\centering
            \footnotesize
            $^e$Centro de investigación biomédica en red\\
            Instituto de Salud Carlos III (CIBER--BBN)\\
            València 46022, Spain
        \end{minipage}
        \hfill
        \begin{minipage}{.45\linewidth}\centering
            \footnotesize
            $^f$Applied Sciences Department,\\
            Northumbria University\\
            Newcastle Upon Tyne NE1 8ST, UK
        \end{minipage}
    \end{minipage}
\end{center}

\medskip
\thispagestyle{empty}

\begin{abstract}
Releasing capsules are widely employed in biomedical applications as smart carriers of therapeutic agents, including drugs and bioactive compounds. 
Such delivery vehicles typically consist of a loaded core, enclosed by one or multiple concentric coating strata. 
In this work, we extend over existing mechanistic models to account for such multi-strata structures, and we characterise the release kinetics of the active substance into the surrounding medium.

We present a computational study of drug release from a multi-stratum spherical microcapsule, modelled through a non-linear diffusion equation incorporating radial anisotropy and space- and time-discontinuous coefficients.
The problem is solved numerically using a finite volume scheme on a grid with adaptive spatial and temporal resolution.
Analytical expressions for concentration and cumulative release are derived for all strata, enabling the exploration of parameter sensitivity --such as coating permeability and internal diffusivity-- on the overall release profile.

The resulting release curves provide mechanistic insight into the transport processes and offer design criteria for achieving controlled release.
Model predictions are benchmarked against \invitro experimental data obtained under physiologically relevant conditions, showing good agreement and validating the key features of the model.
The proposed model thus serves as a practical tool for predicting the behaviour of composite coated particles, supporting performance evaluation and the rational design of next-generation drug delivery systems with reduced experimental effort.

\medskip


\end{abstract}

    \vspace{-1em}
    \tableofcontents
    \vspace{-3.5em}

\newpage

\section{Introduction}
Polymeric microparticles (MPs) are receiving increasing interest as drug delivery systems~\cite{bhu,naz}.
They present interesting features such as different administration routes or encapsulate various bioactive molecules, including antibiotics, chemotherapeutics and nutraceuticals~\cite{birk, wei,dur,puri}.
Particularly, nutraceutical components are credited with a number of health benefits and are associated with market-driven research and product innovation, thus enhancing the quality and safety in the food and health sectors.
However, during production, storage, and even during digestion, these products may be degraded prior to exerting a beneficial effect.
Also, they are often sensitive to light, oxygen and temperature, thus decreasing their shelf life.
In addition, many nutraceuticals are sensitive to the digestive environment and processes, which can affect their bioactivity.
Therefore, understanding the transit, absorption, and bioavailability through the gastrointestinal tract is a keystone for designing tailored releasing MPs for nutraceuticals~\cite{puri}. 

A possible approach to enhance nutraceuticals survivability is given by encapsulation, which preserves and protects them through digestion until they reach the absorption region, where MPs should break down and release their content. 
Currently, smart versatile particles, based on stimuli-responsive materials capable of responding to the biochemical alterations of the environment, such as pH, have attracted great interest due to their structural and morphological advantages. 
Furthermore, the release of the MP payload can be easily tuned by selecting a suitable polymer and its chemical characteristics, such as monomer composition and molecular weight~\cite{lag}.
Another important aspect is the selection of the manufacturing method, according to the drug and specific applications, by which active molecules can be released at a programmed rate to the site of pharmacological action~\cite{gent, desmond2024}.

A matrix platform, where an active agent is uniformly dispersed, is sometimes used.
However, multi-stratum coated systems are nowadays more commonly employed.
These MPs are constituted by a drug-loaded spherical core, surrounded by one or few polymeric concentric strata (shells)~\cite{kaoui}.
In a typical microcapsule (MC), the core has a size of hundreds of micrometres, and the multi-stratum enveloping structure can vary from hundreds of nanometres to a few micrometres.
Such multi-stratum encapsulation enhances the structural stability and biocompatibility, allows the storage of different active agents, and protects the therapeutic agents from external chemical aggression, mechanical erosion, and premature degradation~\cite{gent}. 

Different technologies have emerged in the last years to design and build multi-stratum MCs~\cite{bar,ton,li2017}. 
Various biomacromolecules -- such as proteins, carbohydrates, and their combinations -- have been explored for nutraceutical ingredients protection and stability~\cite{azizi2021}.
Polysaccharides like alginate, chitosan, starch, and carboxymethyl cellulose offer strong film-forming and chemical stability, improving resistance under harsh conditions~\cite{liu2020}.
In particular, alginate is a polysaccharide composed of $\beta$-D-mannuronic acid (M) and $\alpha$-L-guluronic acid (G) residues arranged in varying sequences.
It forms hydrogels in the presence of divalent cations such as Ca$^{2,+}$ through ionic cross-linking, particularly with G-blocks.
Also, alginate is biocompatible, biodegradable, and readily degrades in the intestines, making it suitable for biomedical and food applications~\cite{lee2012}.

While classical \invitro and \invivo techniques experimental techniques remain essential, mathematical and computational modelling offers an alternative tool in drug delivery design and provides more insights on the effects of various geometrical and physical conditions that can be useful to drastically reduce the number of experiments and, also, to provide a better design and manufacturing of drug-loaded systems~\cite{sie}.
To this aim, fully mechanistic continuum models are commonly used, being based on physics-based equations that account for parameters having a direct physico-chemical meaning~\cite{kaoui,bar}. 
In this context, we develop a continuum-scale mathematical model to evaluate the transport and release of combined drug and nutraceutical components from spherical shell-coated MPs immersed in a fluid environment with time-varying pH.
In particular, the model accounts for several different components (core and coating by multiple strata of arbitrary thickness) and considers 
(i) the incorporation of nutraceutical components into the core and 
(ii) the diffusion and the retention of a therapeutic agent from a polymeric MC.

Although diffusion equations have been extensively studied in the literature~\cite{cra}, the size contrast between the core and the shell strata of the MC, as well as the progressive erosion the capsule undergoes, poses interesting numerical questions worth investigation.
To get meaningful results, we will employ a finite volume scheme with variable space and time steps, together with a delicate management of the interfaces between strata, to be interpreted as internal boundary conditions. 
Also, a radial anisotropy is introduced in the diffusion equation to better fit experimental data, especially in the first (internal) strata, which often present particular behaviours (pre-coating).

After a preliminary study of the multi-strata impact on the numerical model, we proceed by identifying the relevant parameters over a set of \invitro experiments,  and we perform a sensitivity analysis to understand the influence of different model conditions and configurations on the release behaviour. 
We show that the model is capable of fitting the observed experimental data and, thus, can be used as a predictive tool. 
The results of our model show good agreement with the data reported in the literature and with the carried out experiments.


\section{The proposed model}\label{sec:model}

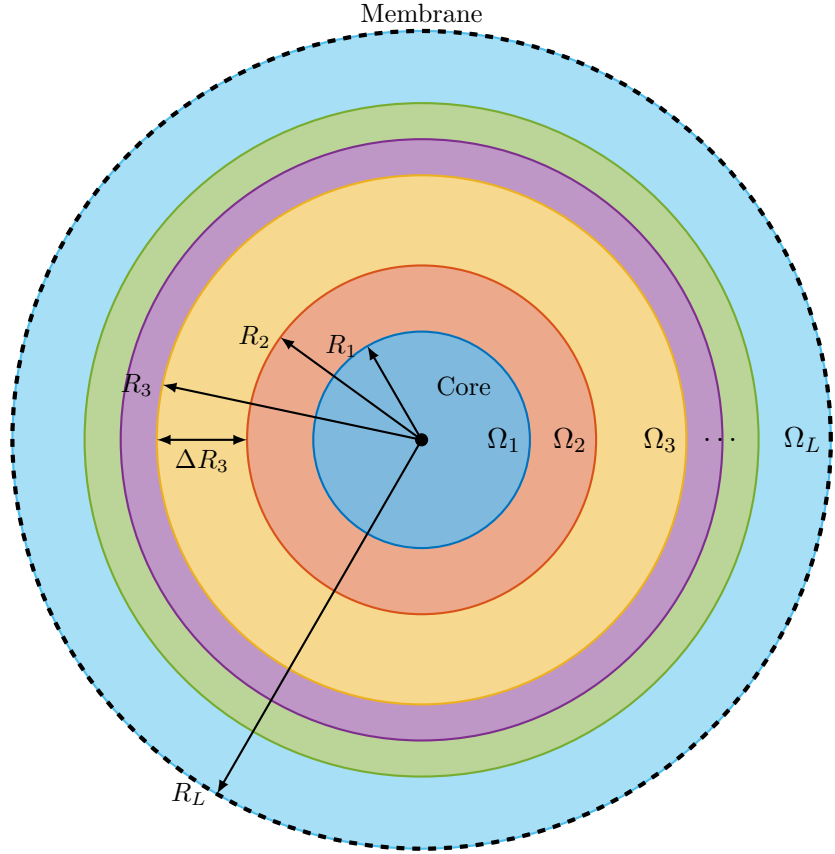
\begin{figure}[ht!]
    \centering
    \input{figure/layerByLayer}
    \caption{
        Cross-section of the multi-stratum MC.
        Figure not to scale.
    }
    \label{fig:LbL}
\end{figure}

Let us consider a multi-stratum spherical microcapsule (MC) made of a core (the depot, $\Omega_1$) surrounded by a sequence of $L-1$ concentric polymeric coating strata ($\Omega_\ell$, with $\ell=2,...,L$), as illustrated in \cref{fig:LbL}.
The core is typically larger in size (hundreds of micrometres) and contains the principal compound like, \eg, drugs or nutraceuticals.
Conversely, the enveloping strata have limited thickness (hundreds of nanometres) and are constituted of different, yet homogeneous and isotropic materials, which increase capsule resistance and structural stability.
In the case of a releasing MC, such strata might be filled with various auxiliary compounds at different concentrations and can be customised to allow selective diffusion, to better control the release rate~\cite{gent}.
The size of these strata may vary considerably from one another, up to some micrometres.
An outer thin protective film, which is generally considered as a membrane, is also included.  
This membrane shields and preserves the encapsulated compounds from degradation and fluid convection, protects the MC overall structure, and guarantees a more controlled and sustained release~\cite{kaoui}.
However, upon ingestion, the membrane immediately dissolves and the outer stratum $\Omega_L$ is considered in contact with the targeted release medium (either a bulk fluid or tissue). During the transit in the digestive tract, the external shells are exposed to the unfavourable and aggressive environment of variable pH and undergo a progressive erosive process.
The biocompound transport within the MC is driven by pure diffusion, but predicting kinetics is not an easy task, due to the multi-stratum structure and the domains thickness contrast.

Furthermore, experiments show that, due to electrostatic interactions with long polymeric chains, a low percentage of travelling molecules remain entrapped and are retained within the polymeric multi-strata without any release.  
We model this aspect by assuming that, while diffusing, part of the drug through the coating can potentially be permanently bound at a given rate~\cite{ton}. 
In other cases, a natural decay of the drug can occur.

\subsection{The general formulation}

Let us consider a general structure made of $L-1$ strata $\Omega_2, \dots, \Omega_L$ surrounding a core $\Omega_1$, hence forming, as a whole, a capsule with a protective shell.
Each domain individually, or a group of domains, is filled with $S$ diffusing substances -- typically, different substances are initially deployed in the core and the strata, yet they can move freely upon ingestion.
Let us denote by $\cs(x,t)$ the concentration of the substance (or phase) $s$ at point $x$ and time $t$.
Note that the position $x$ univocally identifies either a capsule's stratum $\Omega_\ell$, $\ell = 2, \dots, L$, or the core $\Omega_1$.

In the following multi-stratum general model, the biocompound transport in each domain $\Omega_\ell$, $\ell=1,\ldots,L$, can be described by the second Fick's diffusion law and the decay/binding processes as a first-order reaction kinetics:
\begin{equation}
    \frac{\p \cs}{\p t} = \nabla \cdot (D_\ell^s \nabla \cs)- \beta_\ell \Kls
    \qquad
    \mbox{in} \ \Omega_\ell\,,
    \quad
        s=1,\dots,S,
    \label{eq31} 
\end{equation}
where $D_\ell^s$ is the diffusion coefficient and $\beta_\ell^s$ represents the reaction rate of the $s$-th phase in the $\ell$-th domain $\Omega_\ell$. 
The initial concentrations
\begin{equation}\label{def:hatK}
    \cs(\cdot,0)=\csinit(\cdot)   \quad \mbox{in} \ \Omega_\ell 
\end{equation}
are given in each domain, for all substances.

Contiguous domains are paired by imposing flux continuity at the interface, namely:
\begin{equation}
    -D_{\ell-1}^{s} \nabla \cs(x^-) \cdot \hat{n} =  
    -D_{\ell}^{s}   \nabla \cs(x^+) \cdot \hat{n} \qquad  
    \text{ at } x\in\Omega_{\ell-1}\ \cap\ \p \Omega_{\ell}, 
    \qquad \ell=2,\dots,L
    \label{gh3} 
\end{equation}
where $\hat{n}$ is the surface external normal unit vector and $x^\pm$ denotes the inner/outer limit as usual.
Finally, Robin boundary condition is used on the outer stratum, \ie
\begin{equation}
\label{robin}
    \nabla  \cs \cdot \hat{n} = -\lambda^s \cs 
    \qquad  \text{ at } \|x\|=R_L,  
\end{equation}
where 
$\lambda^s \geq 0$  (ms$^{-1}$) is a mass transfer coefficient for substance $s$, expressing the surface coating finite resistance. 
As a special case, $\lambda^s = 0$ (resp., $\lambda^s \to \infty$) represents an impermeable (resp., perfect sink) boundary.

\subsection{The radial symmetric equations}

\begin{figure}[t]
    \centering
\input{figure/layerRadialConstruction.tex}
    \caption{1D radial domain with stratum division for a single substance. Diffusion coefficients depend both on the stratum and the direction of the flux $F$ of the moving substance. The red zone on the far right denotes the leftward erosion, see \Cref{sec:erosion}. 
    }
    \label{fig:setting_analitico_anisotropia_erosione}
\end{figure}
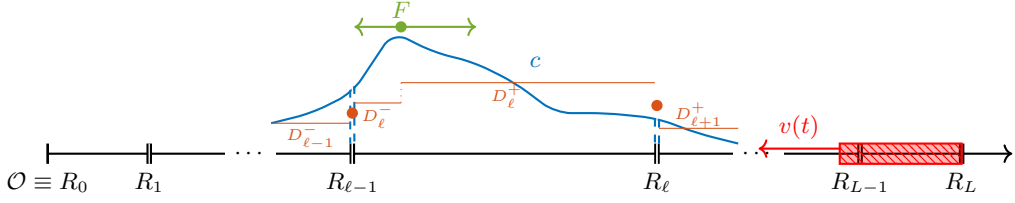

In the above configuration, we assume each stratum has the shape of a concentric spherical shell (annulus) of increasing radius.
Due to the homogeneity, we can assume that net drug diffusion occurs along the radial direction only, and thus we restrict our study to a one-dimensional model for the function $\cs(r,t)$, where $r$ represents the radial coordinate.
Here we denote by $R_\ell$ and $\Delta R_\ell:=R_\ell-R_{\ell-1}$ respectively the radius and the thickness of the of $\Omega_\ell$ (cf.\ \cref{fig:LbL}), where we assume $R_0 = 0$.
Thus, the general formulation given in the previous section reduces to a system of $L \times S$ equations, which in 1D radial symmetry reads, for any $\ell=1,\dots,L$ and $s=1,\dots,S$, as 
\begin{align}
    &\frac{\p \cs}{\p t} = \frac{1}{r^2} \frac{\p}{\p r}\left(D_\ell^s\left(\frac{\p \cs}{\p r}\right)r^2 \frac{\p \cs}{\p r}\right) - \beta_\ell^{s} \cs\ ,
    &\begin{array}{l}
        r \in (R_{\ell-1}, R_\ell)\\
    \end{array}
    \label{eq:rs1}
    \\
    &\cs(r, 0) = \csinit(r)\ ,
    &\begin{array}{l}
        r \in (R_{\ell-1}, R_\ell)\\
    \end{array}
    \\
    &\frac{\p \cs}{\p r} = 0\ ,
    &\begin{array}{l}
        r = 0\\
    \end{array}
    \\
    &-D_{\ell-1}^{s}\frac{\p \cs(r^-)}{\p r} = -D_{\ell}^{s}\frac{\p \cs(r^+)}{\p r}\ ,
    &\begin{array}{l}
        r = R_{\ell-1}\\
    \end{array}
    \\
    &\frac{\p \cs}{\p r} = -\lambda \cs\ ,
    &\begin{array}{l}
        r = R_{L}\\
    \end{array}
    \label{eq:rsn}
\end{align}
where $\csinit$ is the initial concentration defined in \eqref{def:hatK}. 
Note that in \eqref{eq:rs1} we have included the dependence of $D$ on the derivative of $c$ in $r$. 
This was done to allow a possible radial anisotropy in the diffusion dynamics, meaning that, for a given $r$, the diffusion towards the centre can be different from the diffusion towards the outer boundary, see also \cite{xu}. 
In the following experiments, we simply define it as a piecewise constant value, depending on the flux direction only, \ie for $\ell = 1, \dots, L$:
\begin{equation}\label{def:D+-}
    D_\ell^s\left(\frac{\p c^s}{\p r}\right) =
    \left\{
    \begin{array}{ll}
        D_\ell^{s,+} & \text{if } \ds{\p c^s \over \p r} \geq 0\\ [2mm]
        D_\ell^{s,-} & \text{if } \ds{\p c^s \over \p r}< 0 
    \end{array}
    \right.
\end{equation}
where we define the anisotropic factor $\alpha_\ell^s \geq 0$ as the ratio $\alpha_\ell^s := D_\ell^{s,-} / D_\ell^{s,+}$ (here $\alpha=0 / +\infty$ would describe a material admitting outward/inward flux only, $\alpha=1$ an isotropic material).
\cref{fig:setting_analitico_anisotropia_erosione} shows the 1D domain with the stratum division and how to manage the radial anisotropy.

\paragraph{Computing the released mass}
To characterise the release process, we evaluate the cumulative mass $m$ of the released compound $s$ as a function of time.
This quantity can be obtained by integrating the corresponding concentration fluxes over the outer surface of the MC.
In formulas, we obtain
\begin{equation}
    m^s(t) = \ds\int_0^t \left( \ds\int_{\p \Omega_L} -D_L^{s} \nabla \cs \cdot \hat{n}  \right) ds \ ,
    \qquad s=1,\dots,S \ ,
    \label{eq:released}
\end{equation}
which, under the radial symmetry, reduces to
\begin{equation}
    m^{s}(t) = - \ds\int_0^t D_L^{s} R_L^2 \frac{\p \cs(R_L,s)}{\p r} ds \ ,
    \qquad s=1,\dots,S \ .
    \label{eq:rs-release}
\end{equation}

\subsection{Modelling the erosion}\label{sec:erosion}

So far, we have assumed the domain of the MC remains unchanged throughout the process.
However, because of water imbibition and chemical interactions with the release medium, or due to the unfavourable pH of the environment, the MC undergoes a degradation/erosion process on its surface.
As a consequence, the external stratum erodes and its radius reduces over time, progressively melting the strata until the core itself becomes exposed. At that point, the whole capsule, left without any protective shell, typically disintegrates.

From the modelling perspective, let $R(t)$ be the radius of the sphere at time $t$.
The speed of erosion can be modelled as some positive function $v_e(t)$ such that
\begin{equation}\label{def:R(t)}
    R(t) := \max\left\{R_1, R_L - \int_0^t v_e(s)s\ ds\right\}, \qquad t \geq 0.
\end{equation}
The simulation domain shrinks progressively, and the domain for \eqref{eq:rs1} becomes $(R_{\ell-1}, \min\{R_\ell, R(t)\})$, see \cref{fig:setting_analitico_anisotropia_erosione}.
As a consequence, when $R(t) \leq R_{\ell-1}$, then $\Omega_\ell$ is considered completely eroded and its contribution is removed from the model;
hence, we denote with $\ell_L$ the most external non-eroded domain.
Note that we will consider the function $v_e(t)$ as assigned or estimated from experiments. 

Clearly, if \eqref{eq:rs1} is solved assuming any time a fixed domain, then an additional term is required to take into account the mass instantaneously released due to volume erosion.
Hence \eqref{eq:rs-release} reads as:
\begin{equation}
    m^{s}(t) = - \ds\int_0^t D_{\ell_L}^{s} R(s)^2 \frac{\p \cs(R(s),s)}{\p r} ds + 
        \int_0^t \cs(R(s),s) 4\pi R^2(s)v_e(s) ds \ ,
    \label{eq:rs-release-complete}
\end{equation}
where, for the additional term, we have considered that the infinitesimal mass element $dm$ is given by the concentration $c^s$ multiplied by the infinitesimal annulus volume element swept by $R$ in time $dt$,  $dV =\frac43 \pi ((R+dR)^3-R^3)\approx 4\pi R^2 dR$ (first order approximation), and, finally, $dR=v_e dt$.

\section{The numerical model}

In this section, we detail the numerical method used for approximating the system \eqref{eq:rs1}--\eqref{eq:rsn}.
Due to strata thickness contrast (several orders of magnitude between the core and external strata), a numerical grid with variable space and time steps is required. 
Moreover, the erosion moves the MC surface inwards and  requires attention to the right boundary condition.

Since the kinetics of concentrations are independent of one another, and all undergo similar evolution in space and time, in what follows, we will refer to a generic concentration $c$, omitting the $s$ index.

\subsection{The numerical scheme for a fixed domain}\label{sec:computationa-model-fix}

We preliminarily assume that we are dealing with a \textit{time-independent} spherical domain, \ie where $R(t)=\bar R = R_L$, $\forall t$.
We assume $t_0=0$ as the initial time and $r=0$ is the centre of the \insilico capsule. 
The numerical solution is found in the spatio-temporal domain $[0,\bar R] \times [0,T]$, by tracking the evolution of $c(r,t)$ and the corresponding released mass $m(t)$.
For numerical convenience, we divide each stratum $\ell$ into a number of grid cells, say $N_\ell$.
In each stratum, the length of cells is constant and given by $\Delta r^\ell$.
The centre of the $j$-th cell of the $\ell$-th stratum is denoted by $r_j^\ell$ and its value corresponds to the distance from the centre of the capsule to the centre of the cell itself.
Following the same notation, we denote by $\big[r^\ell_{j-\sfrac12}, r^\ell_{j+\sfrac12}\big)$ the $j$-th cell of the $\ell$-th stratum, with $r^\ell_{j\pm\sfrac12}:=r_j^\ell\pm\frac{\Delta r^\ell}{2}$.
Finally, we denote with $\dt$ the time step, and its number with $n$, as usual.

The numerical solution will be given by the values $C_j^{\ell,n}$, for all $\ell$, $j$, and $n$, which approximate the average concentration in each cell at each time step in each stratum \cite{levequebook},

\be
    C_j^{\ell,n} \approx \frac{1}{\Delta r^\ell}\int_{r^\ell_{j-\sfrac12}}^{r^\ell_{j+\sfrac12}} c(r, n\, \dt) dr. \label{eq4}
\ee

We recall that, due to the spherical symmetry, the interfaces between cells are spherical shells of surface
\be
    A_{j+\sfrac12} = 4 \pi r_{j+\sfrac12}^2,  \label{eq5}
\ee
and that the $j$-th cell accounts for a total volume of 
\be
    V_j = \frac43\pi(r_{j+\sfrac12}^3-r_{j-\sfrac12}^3). \label{eq3}
\ee

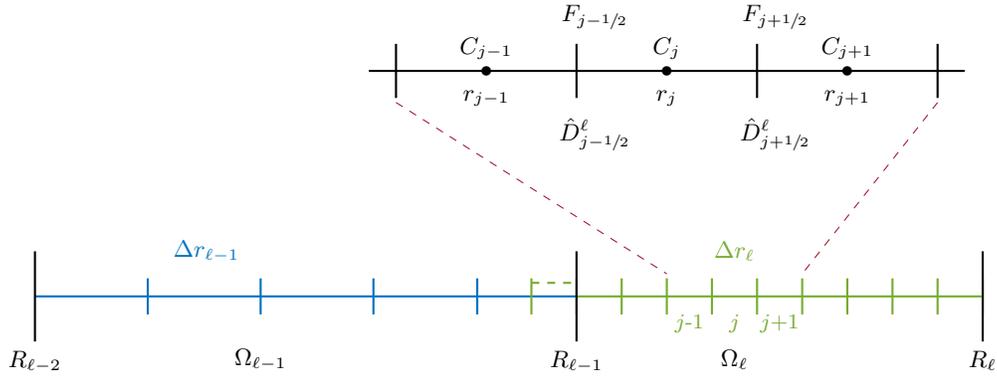
\begin{figure}[h!]
    \centering
    \input{figure/numericScheme}
    \caption{
    Grids and variable settings. Coarse grid on $\Omega_{\ell-1}$ and fine grid on $\Omega_{\ell}$. Inner cells zoom: diffusion coefficients $\hat D$ (\cref{def:hatD}) and numerical fluxes $F$ (\cref{numerics-diffusion}) are defined at the cell interfaces, diffusion coefficients $D$ and concentrations $C$, instead, are defined at the cell centres. The dashed cell is the ghost cell for the fine grid.
    } 
    \label{fig:domain-numerical-setting}
\end{figure}

The stratum-dependent diffusion coefficients can be discretised at interfaces using the harmonic mean,
\begin{equation}\label{def:hatD}
    \hat D^\ell_{j+\sfrac12} := 
    \begin{cases}
        D_\ell^+, & \text{if } C_j^\ell \geq C_{j+1}^\ell \text{ and } R_{\ell-1} \leq r_{j+\sfrac12} < R_\ell,\\
        D_\ell^-, & \text{if } C_j^\ell < C_{j+1}^\ell \text{ and } R_{\ell-1} \leq r_{j+\sfrac12} < R_\ell,\\
        \ds\frac{2D^+_\ell D^+_{\ell+1}}{D^+_\ell + D^+_{\ell+1}}, & \text{if } C_j^\ell \geq C_{j+1}^\ell \text{ and } r_{i+\sfrac12} = R_\ell,  \\
        \ds\frac{2D^-_\ell D^-_{\ell+1}}{D^-_\ell + D^-_{\ell+1}}, & \text{if } C_j^\ell < C_{j+1}^\ell \text{ and } r_{i+\sfrac12} = R_\ell.
    \end{cases}
\end{equation}
From now on, when possible, we drop the index $\ell$ for readability.
In \cref{fig:domain-numerical-setting}, we depict the numerical setting at the junction between a sample coarse and fine grid, and we visualise the quantity $\hat D$ defined at interfaces.

\paragraph*{Diffusion in spherical symmetry}

The second-order operator in \eqref{eq:rs1} can be approximated by means of the following finite volume discretisation
\be \label{numerics-diffusion}
    \begin{split}    
        &\frac1{r^2} \frac{\partial}{\partial r} \left( D\ r^2 \frac{\partial c}{\partial r}  \middle)\right|_{r=r_j}\approx
        \\
        &\qquad\approx\frac{1}{r_j^2}
        \frac{1}{\dr} \left(
            \underbrace{\hat D_{j+\sfrac12}\ r_{j+\sfrac{1}{2}}^2\frac{C_{j+1}^n - C_j^n}{\dr}}_{F_{j+\sfrac12}}
            -
            \underbrace{\hat D_{j-\sfrac12}\ r_{j-\sfrac{1}{2}}^2\frac{C_j^n - C_{j-1}^n}{\dr}}_{F_{j-\sfrac12}}
        \right)
    \end{split}
\ee
where $F_{j+\sfrac12}$ and $F_{j-\sfrac12}$ are the right-outward and left-inward numerical fluxes of the cell $j$, respectively (cf.\ \cref{fig:domain-numerical-setting}).

\paragraph*{Time derivative}
For what concerns the time derivative $\ds\frac{\partial c}{\partial t}$, a natural choice is forward finite differences, 
\be
    \left.\frac{\partial c}{\partial t}\right|_{r=r_j} \approx 
    \frac{C_j^{n+1} - C_j^n}{\Delta t}.
\ee

\paragraph*{Complete scheme in each stratum}
Overall, the complete numerical scheme reads
\be\label{eq:complete scheme}
        C_j^{n+1} = C_j^n + \frac{\dt}{\dr \ r_i^2}(F_{j+\sfrac12}-F_{j-\sfrac12}) - \dt\beta_j C_j^n, 
\ee
where the last term accounts for the decay process and it is considered fully independent of space and of the diffusion process.
The above numerical scheme is stable only under the CFL condition 
\cite{levequebook}:
\be
    \dt \leq \frac{\dr^2}{2\max_j(\hat D_{j\pm\sfrac12})}.
\ee

\paragraph*{Stratum interfaces and variable space-time steps}
Some important remarks are in order:
\begin{itemize}
    \item The big difference between the thickness of the core and that of the coating strata suggests using different space steps.
To this end, it is convenient to assume that space steps $\Delta r$ are constant in each stratum and are all multiples of the smallest. 
This makes computation easier to manage.
Moreover, when two strata with different space steps join together, we compute the mass balance at the interface only on the side of the finer grid. 
To do that, a ghost cell, with a size equal to the small cell, is placed on the side of the stratum with large cells (cf. \cref{fig:domain-numerical-setting}).
The value of the concentration of the ghost cell is assumed to be equal to the value of the large cell in which it is embedded.
Once the mass balance in the small cell is computed, it is used (with opposite sign) for evaluating the mass balance in the large cell.
\item The CFL condition \cite{levequebook} must be kept as sharp as possible (ideally it should be satisfied with equality), in order to avoid numerical diffusion (for computational accuracy) and to avoid small time steps (for computational efficiency).
Therefore, strata with different space steps must have different time steps, which means that not all strata are updated with the same frequency.
Coming back to the mass balance at stratum interfaces discussed before, different time steps require storing in a buffer the mass which must pass from the fine to the coarse layer, in such a way that the right amount of mass is transferred when the coarse grid is eventually updated.
\item The large difference between space steps introduces some unfeasible dynamics in the numerical approximation. 
For example, when some mass is transferred from the fine to the coarse grid, the mass spreads instantaneously along the whole large cell (since at each given time the concentration value is constant in each cell), thus covering a distance which is too large to be consistent with the physics of the problem. 
On the other side, when some mass is transferred from the coarse to the fine grid in one time step (the large one), the mass has enough time to travel through the entire neighbouring small cell and reach the next small cell, but this is not allowed by the CFL condition, being a 3-point stencil used here (fig. 3). 
To overcome this problem, we introduce one or more domains (fictitious strata), which allows us to move \textit{continuously} from the coarse to the fine grid (or vice-versa), thus avoiding excessive differences in the grid size.
The model, in fact, includes the possibility of having a large number of domains, not only to describe the experimental physical strata with due accuracy, but also for numerical reasons.
\end{itemize}

\paragraph*{Boundary conditions at $r=\bar R$} 
At the outer boundary, Robin condition (see \cref{robin}) is imposed by means of an additional ghost cell, indexed by $N_L+1$, where we set
\begin{equation}\label{robinBC}
    C_{N_L+1}^n=C^n_{N_L}\left(1-\Delta r^L\lambda\right) \ ,
\end{equation}
corresponding to an outgoing flux (cf.\ \cref{numerics-diffusion}) through the last interface $N_L+\sfrac12$ given by
\begin{equation}
    F_{N_L+\sfrac12} = \hat D_{N_L+\sfrac12}\ r^2_{N_L+\sfrac12}\ \lambda C_{N_L}\ .
\end{equation}

\paragraph*{Modelling the erosion}\label{sec:numerics-erosion}

So far, we have discussed the problem by considering a fixed spatial domain.
Now we generalise the model accounting for the erosion term, as defined in \eqref{def:R(t)}.
In order to keep the numerical implementation as simple as possible, we have assumed that the external boundary condition \eqref{robinBC} progressively moves inward.
Moreover, we never allow a cell to be only partially eroded. 
Instead, we reduce the domain by one cell only when it is fully disintegrated.
When this happens, all the mass contained in that cell is released, and the ghost cell accounting for the external boundary condition is shifted one cell inwards.

\section{Numerical validation }

In this section, we introduce a few simple tests to prove the model effectiveness and to validate the numerical method. 
Let us consider a homogeneous sphere of $R = 100$ $\mu$m, filled with a given substance with constant density of $c_0 = 1$ $\mu$g/$\mu$m$^3$, and let us assume a diffusion coefficient $D = 0.01$ $\mu$m$^2$/s.
For the sake of simplicity, let us assume the substance being non-decaying ($\beta = 0$ s$^{-1}$) and the MC is not eroding ($v = 0$ $\mu$m/s) and isotropic ($\alpha = 1$).
The total mass ($m = 4.188$ g) should then completely leave the medium (up to $\sim99.9\%$) within 4 hours ($T=14400$ s). Parameters are summarised in \cref{tab:parameters-preliminary-tests}.
\begin{table}[h!]
    \setlength{\tabcolsep}{8pt}
    \centering
    \scalebox{.8}{
    \begin{tabular}{|c|c|c|c|c|c|c|c|}
        \hline
        $R$ [$\mu$m]
        & $T$ [min]
        & $c_0$ [$\mu$g/$\mu$m$^3$]
        & $D$ [$\mu$m$^2$/s]
        & $\beta$ [s$^{-1}$]
        & $v$ [$\mu$m/s]
        & $\lambda$ [$\mu$m$^2$/s ]
        & $\alpha$ [-] \\\hline
        100 
        & 240 
        & 1 
        & 0.5 
        & 0 
        & 0 
        & $1$
        & 1 \\\hline
    \end{tabular}
    }
    \caption{Set of parameters used in preliminary tests}
    \label{tab:parameters-preliminary-tests}
\end{table}

In order to capture possible numerical differences between mono- and multi-stratum systems, we run the same physical model under a variety of different grid combinations:
\begin{description}
    \item[one-stratum, extremely fine grid, reference case.] We set a single extremely fine grid as a reference (exact) solution, setting $\dr = 0.01$ $\mu$m and setting the corresponding time discretization to $\dt = 10^{-5}$ s, a reasonable value according to the CFL condition required for the stability of the scheme.
    \item[one-stratum, fine grid.] We set a single fine grid as the best case benchmark, setting $\dr=0.1$ $\mu$m and $\dt=10^{-3}$ s.
    \item[one-stratum, coarse grid.] We set a single coarse grid as a worst-case benchmark, setting $\dr=1$ $\mu$m and $\dt=0.1$ s.
    \item[10-strata, constant coarse grid.] We run the simulation over 10 evenly-spaced strata, employing the very same discretisation as in the coarse grid example for each stratum, so as to verify the impact of the multi-grid over-engineering.
    \item[10-strata, variable-time fine grid.] As before, we adopted 10 strata and we fixed the spatial discretisation to the one of the fine grid, \ie $\dr=0.1$ $\mu$m; however, we set variable $\dt$s, ranging between $5\cdot 10^{-5}$ and $10^{-3}$ to understand the impact of different temporal discretisation on the results. In particular, we set both progressive refinements and coarsening of $\dt$ between strata as follows:
    $$
        \dt = 10^{-3}\cdot [
            1,   0.5,    0.1,    0.05,    0.1,    0.5,   1,    0.05,    0.05,    1   ]
    $$
    \item[2-strata, coarse-fine coupled grid.] We coupled the coarse grid and the fine grid by running the simulation over an internal stratum of 75 $\mu$m coarse-discretised and an external stratum of 25 $\mu$m fine-discretised.
    \item[2-strata, fine-coarse coupled grid.] By reversing the previous configuration, we finally cou\-pled an external coarse grid with an internal fine grid.
\end{description}

\begin{figure}[t]
    \centering
    \includegraphics[width=0.9\linewidth]{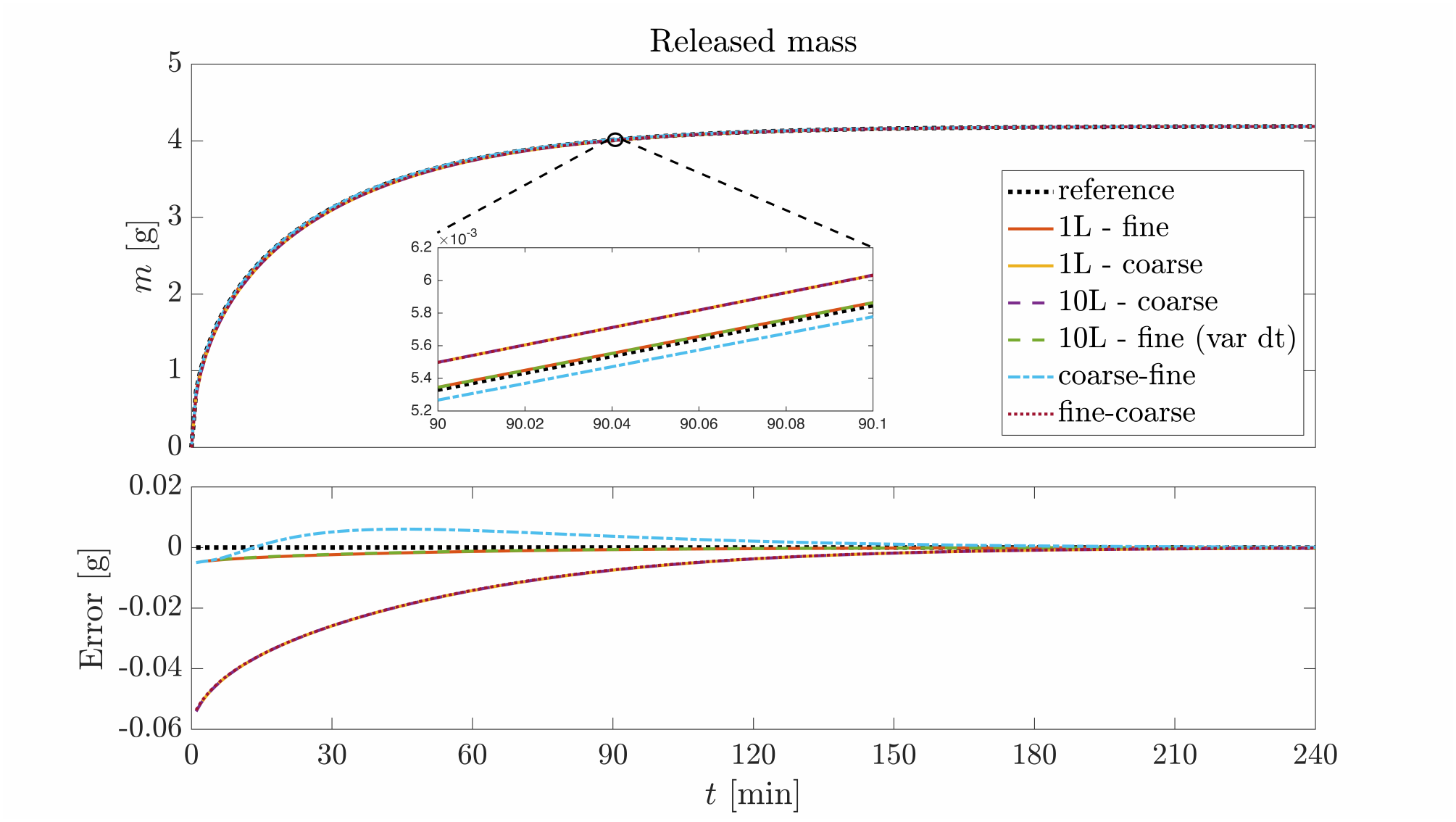}
    \caption{Released cumulative mass (top) and corresponding signed absolute error (bottom) as a function of time, for the reference solution and the other 6 configurations.}
    \label{fig:academic-release} 
\end{figure}
\begin{figure}[t]
    \centering
    \includegraphics[width=0.7\linewidth]{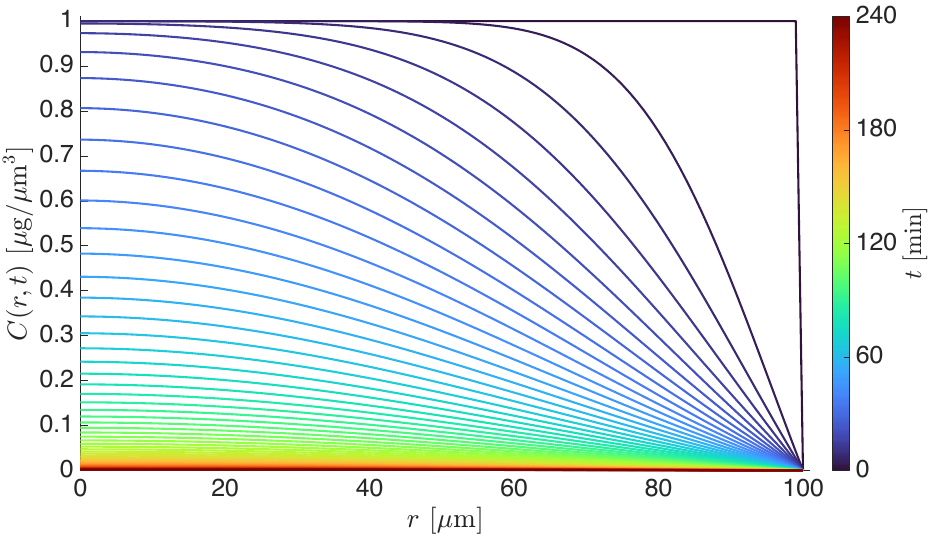}
    \caption{
    Concentration profiles at equidistant time intervals (240 s)
    along the radius. The colour bar indicates the time scale.
    }
    \label{fig:academic-profile}
\end{figure}

\begin{table}[t]
    \centering
    \begin{tabular}{c|c|c|c||c}
        \textbf{Number of strata} & \textbf{Grid} & $\dr$ [$\mu$m] & $\dt$ [s] & Mean relative error [$\%$]\\
        \hline
        1  & fine        & 0.1   & 0.001    & $0.034$ \\
        1  & coarse      & 1     & 0.1      & $0.376$  \\
        10 & coarse        & 1   & 0.1      & $0.376$ \\
        10 & fine (var $\dt$)    & 0.001    & 0.005 -- 0.1    & $0.033$ \\
        2  & coarse+fine & 1/0.1 & 0.02/0.001 & $0.073$ \\
        2  & fine+coarse & 0.1/1 & 0.001/0.02 & $0.374$ \\
    \end{tabular}
    \caption{Relative error between the six different configurations compared with the reference case. In all cases, the error is relatively small.}
    \label{tab:academic-release}
\end{table}

In \cref{fig:academic-release}, we plot the release curve (top) and the absolute error (bottom, point-wise signed distance from the reference), for released mass computed every minute.
Due to the absence of the binding effect, the whole initial mass is released within 240 minutes, as expected, also explaining why the error goes to zero in the limit, for each scenario.
In \cref{tab:academic-release} , we report the corresponding discrepancies (vs.\ the reference mono-stratum finer grid approximation of the exact solution) along with the experiments summary.
\cref{fig:academic-profile} shows the concentration profiles along the radius at time intervals of 4 minutes. 
As expected, the concentration decays to zero after $\sim4$ h, with the release rate larger at early instants.

Unsurprisingly, the fine-grid solution minimises discrepancy over the other experiments.
However, the impact of the sole multi-stratum is negligible w.r.t.\ the single-stratum experiment (cf.\ coarse vs.\ coarse multigrid).
Furthermore, the gain obtained from the improvement in $\dt$ frequency in a multi-stratum environment actually overcomes the performance degradation introduced by the multiple strata themselves (cf.\ fine grid vs.\ multigrid fine).
Finally, the disposition of the strata impacts the simulation accuracy, mainly depending on the flux direction (\ie mass moving from coarse to fine or vice-versa).

\section{A case study: the coated microcapsules}

In this section, we first describe the experimental setup of the case study involving coated nano-engineered microcapsules, which serves as the basis for the numerical model introduced earlier.
We then demonstrate the ability of the proposed model to reproduce the experimental data with good agreement.
Finally, we present a sensitivity analysis to explore the influence of key model parameters.

\subsection{Experimental setup}\label{sec:experimentalsetup}
To ground our modelling approach in a realistic application, we consider a case study involving microcapsules nano-engineered via multi-strata coating.
In what follows, we detail the experimental procedure adopted for the fabrication and characterisation of these microcapsules, which provides both the physical context and quantitative data used to inform and validate our numerical model (statistics are presented as average and variance over 9 experimental repetitions).

\begin{figure}[t]
    \centering
    \includegraphics[width=0.95\linewidth]{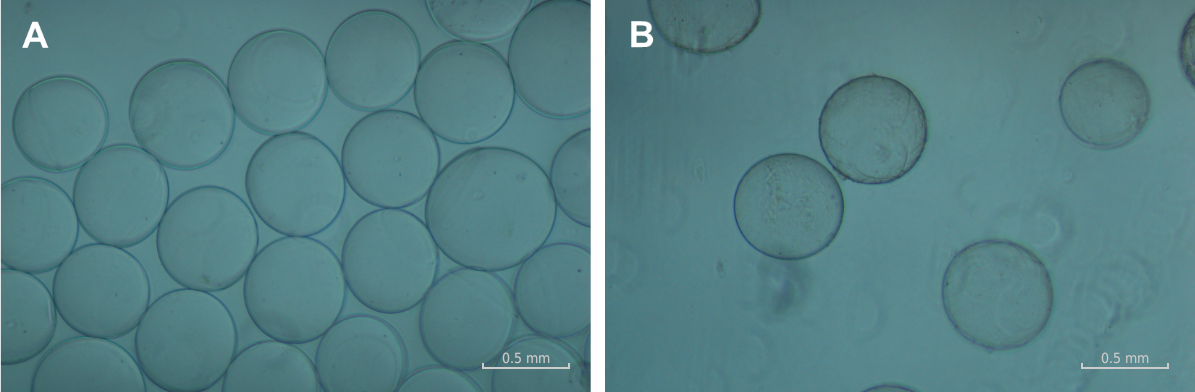}
    \caption{Optical images of MCs: without (A) and with (B) the polysaccharide-based multi-strata shell.}
    \label{fig:beadsAB}
\end{figure}

For the core formation, alginate powder was dissolved in deionised water to prepare a 2 wt\% solution.
Then, the solution was transferred into a syringe and extruded through a 27G conical tip into a 1\% (w/v) CaCl$_2$ gelling solution.
To achieve MC, an adapted electrospinning setup (Spinbox, Bioinicia, Valencia, Spain) was employed. 
The process parameters included a flow rate of 500 $\mu$L/h, an applied voltage of 6.5 kV, and a 2 cm distance between the tip and the gelling bath.
The MC diameter and circularity were evaluated using microscopy (Nikon Eclipse TS100), which revealed an average radius of 285 $\pm$ 18 $\mu$m and a circularity of 0.98 $\pm$ 0.02, see \cref{fig:beadsAB}A. 
Capsules were subsequently coated via coupled deposition of polysaccharides embedding the tannic acid (TA; Merck, UK) as selected therapeutic agent.
TEM (Philips CM 100 Compustage FEI) was used at a voltage of 100 kV for investigating the thickness of each stratum produced, revealing an average value of 18 nm.

During the initial deposition cycles, a couple of pre-coating strata typically forms, characterised by irregular adsorption and partial surface coverage.
This early phase is crucial for modifying the surface chemistry and charge of the substrate, thereby improving its compatibility with the subsequent growth of ordered strata.
As the deposition progresses, the pre-coating serves as an intermediate interface that promotes more uniform molecular organisation, eventually leading to the formation of a continuous, homogeneous multi-strata film with well-defined thickness~\cite{Decher92, Porcel07}.
This transition from a disordered or patchy initial stratum to an ordered stratified structure is particularly evident when working with low-reactivity or rough substrates.

Finally, the MCs were frozen in liquid nitrogen and stored at -80 $^\circ$C for further analysis.
No significant changes were observed in the size of the capsules due to the nanometric thickness of the resulting shell, see \cref{fig:beadsAB}B.

To estimate the content of TA encapsulated within the shell, a dissolution test was conducted in PBS (pH 7.4, 37 $^\circ$C) for 48 hours to allow the degradation of the multi-stratum. 
The released TA solution was evaluated spectrophotometrically at 280 nm (FLUOstar Omega, BMG Labtech, Germany) and revealed a content of 4.721 $\pm$ 0.67 mg/mL in the coated MCs, corresponding to a mass of $m^{\text{TA}}=472.10 \pm 67$ $\mu$g.

Then, \invitro digestion experiments were carried out according to the standardised INFOGEST 2.0 method, comprising sequential oral, gastric, and intestinal phases under physiologically relevant conditions~\cite{brodkorb2019}. 
In the oral phase, samples were mixed with simulated salivary fluid (SSF) containing $\alpha$-amylase and incubated at 37 $^\circ$C for 2 min. 
The gastric phase was performed using simulated gastric fluid (SGF) with pepsin (pH 3.0) and incubated at 37 $^\circ$C for 2 h in a TSSWB15 shaking water bath (ThermoScientific, UK). 
The intestinal phase followed, using simulated intestinal fluid (SIF) containing pancreatin and bile salts (pH 7.0), also incubated at 37 $^\circ$C for 2 h with constant shaking. 
Enzyme activities and electrolyte concentrations were adjusted according to INFOGEST guidelines.

Aliquots (50 $\mu$L) were withdrawn from the digestion mixtures at each phase (oral, gastric, and intestinal), specifically every 30 min during the gastric phase and every 15 min during the intestinal phase, to monitor the release of tannic acid from the beads (see also later Figure~\ref{fig:experiment-release}).
Negligible release was detected during the oral phase, indicating the multi-stratum coating effectively protected the core under salivary conditions.
During the gastric phase, $\sim21.92\%$ of the TA was released, likely due to partial disruption of the coating in the acidic environment.
A more sustained release was observed in the intestinal phase, attributed to the combined effects of hydrolytic degradation and diffusion through the multi-stratum structure.

\begin{figure}[t]
    \centering
    \includegraphics[width=0.9\linewidth]{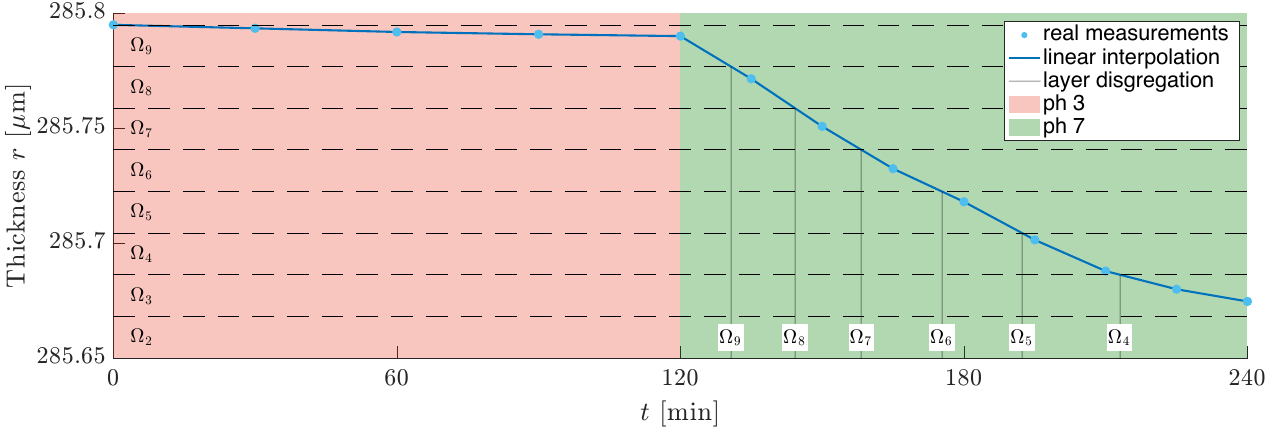}
    \caption{Radius of the eroding MC as a function of time. Radius shrinks due to the erosion of the strata of the shell, which, in turn, depends on the pH during digestive transit.}
    \label{fig:experiment-disgregation}
\end{figure}

Concurrently with aliquot sampling, TEM nanometric measurements on the multi-stratum thickness were carried out to characterise the impact of the erosion process.
Negligible erosion was detected ($<0.5\%$) during the oral phase, further confirming the coating effectiveness.
Interestingly, extremely constrained erosion, corresponding to $3.24\%$ of the total thickness, was measured during the gastric phase, prompting the release of the TA being ascribable to diffusive processes mainly.
Conversely, the consistent erosion during the intestinal phase positively correlates with the higher measured TA releases, suggesting it as the main contributor in the release process.
Also noteworthy is the sudden change within erosion speed upon reaching the last two strata, see \cref{fig:experiment-disgregation}, further suggesting a structural difference between the pre-coating and the actual MC structure.

These results suggest that the encapsulation system provides controlled release, with minimal premature leakage and enhanced delivery in the target intestinal environment.

\subsection{Numerical results}
\label{ssec:estimating-parameters}

The MC model presented in \cref{sec:model} contains a number of interdependent parameters, defined within different subdomains, most of which are not known a priori.
Moreover, even when estimates of certain parameters are available from the literature, these are often subject to high variability and uncertainty.
As a consequence, obtaining reliable estimates is a significant challenge while modelling the release of MCs. 
Starting from a wide range of physically realistic parameters, we then address this issue in two steps:
firstly, by comparison with experimental data sets, we inversely estimate the parameters for the specific MC release system considered and described in \cref{sec:experimentalsetup};
secondly, we studied the parameters sensitivity around the estimated reference values to gain insights on how each component influence the overall release phenomenum.

More in details, using the experimental values for TA release and the erosion curve from \cref{fig:experiment-disgregation}, we identify the optimal design parameters $D^+, \alpha, \lambda, \beta$ for each stratum that best fit the experimental dataset.
In particular, due to the negligible effect of the oral phase (2 min), we simulated the evolution within gastric (2 h) and intestinal (2 h) phases only, for a total of $T=14,\!400$ s.
The reference parameters we obtained are reported in \cref{tab:parameters-real-test}.

\begin{table}[t]
    \centering
    \setlength{\tabcolsep}{8pt}
    \renewcommand{\arraystretch}{1.3} 
    \begin{tabular}{|l|c|c|c|c|c|}
       \hline
       \multirow{3}{*}{\textbf{Parameter}} & \multicolumn{3}{c|}{\textbf{Core} $\Omega_1$} & \multicolumn{2}{c|}{\textbf{Shell} $\Omega_2$--$\Omega_9$}\\
       & \multicolumn{3}{c|}{\textbf{(numerically partitioned)}} & \multicolumn{2}{c|}{\textbf{(8 strata)}}\\
       & $\Omega_{1}^a$ & $\Omega_{1}^b$ & $\Omega_{1}^c$ & $\Omega_2$--$\Omega_3$ & $\Omega_4$--$\Omega_{9}$ \\\hline
       $R$ [$\mu$m]        & \multicolumn{3}{c|}{$285.65$} & \multicolumn{2}{c|}{$285.794$}
       \\\hline
       $\Delta R$ [$\mu$m] & $280$ & $5$ & $0.65$ & $0.0180$ & $0.0180$
       \\\hline
       $c_0^\text{TA}$ [$\mu$g/$\mu$m$^3$] & \multicolumn{3}{c|}{0} & $5.085 \times 10^{-3}$ & $2.543 \times 10^{-3}$\\\hline\hline
       $\lambda$ [$\mu$m  s$^{-1}$] & \multicolumn{5}{c|}{$0.05$}\\\hline  
       $\beta$ [s$^{-1}$] & \multicolumn{5}{c|}{$0$}\\\hline
       $D^+$ [$\mu$m$^2$/s] & \multicolumn{3}{c|}{$6 \times10^{-7}$} & $5\times10^{-6}$ & $10^{-6}$\\\hline
       $\alpha$ [dimless] & \multicolumn{3}{c|}{$0.5$} & $0.2$ & $1.0$ \\\hline\hline 
       $\dt$ [s] & $1$ & $0.05$ & \multicolumn{3}{c|}{$0.01$}\\\hline
       $\dr$ [$\mu$m] & $35$ & $0.5$ & $0.005$ & \multicolumn{2}{c|}{$0.001$}\\\hline
    \end{tabular}
    \caption{Set of parameters used in the case study. MC's core is divided into 3 strata, only for computational purposes. From top to bottom, the first group of parameters are experimentally measured, the second group are best-fitted parameters matching experimental data, and the third group includes numerical parameters.}
    \label{tab:parameters-real-test}
\end{table}

\begin{figure}[t]
    \centering
    \includegraphics[width=0.99\linewidth]{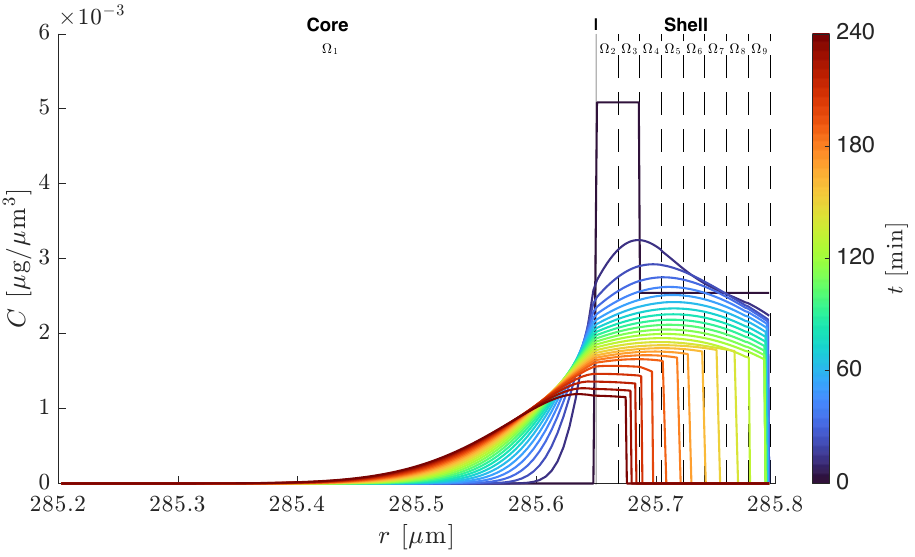}
    \caption{
        TA concentration profiles at equidistant time intervals (10 min) along the radius on the external part of the MC.
        Concentration vertical lines correspond to the progressive stratum breakup in consequence of the erosion (cf.\ vertical lines in \cref{fig:experiment-disgregation}).
        The colour bar indicates the time scale.
    }
    \label{fig:experiment-profile}
\end{figure}

\cref{fig:experiment-profile} shows the TA concentration profiles up to 240 minutes at time intervals of 10 minutes in the outer part of the core and in the strata of the shell. 
The initial difference in concentration between pre-coating ($\Omega_2$ and $\Omega_3$) and the other strata progressively smooths out over time, yet allowing a part of the TA filtrate within the core $\Omega_1$.
Upon erosion (cf.\ progressive vertical drop in concentration), the shell progressively empties, up to when the concentration within the core overcomes the one in the pre-coating;
as a consequence, in the later steps of the experiment, the core serves as a reservoir to keep coating concentration high, and allowing a smoother release when paired with a higher outgoing diffusion (cf.\ \cref{tab:parameters-real-test}, $\alpha_2, \alpha_3 = 0.2$). 

\begin{figure}[t]
    \centering
    \includegraphics[width=0.99\linewidth]{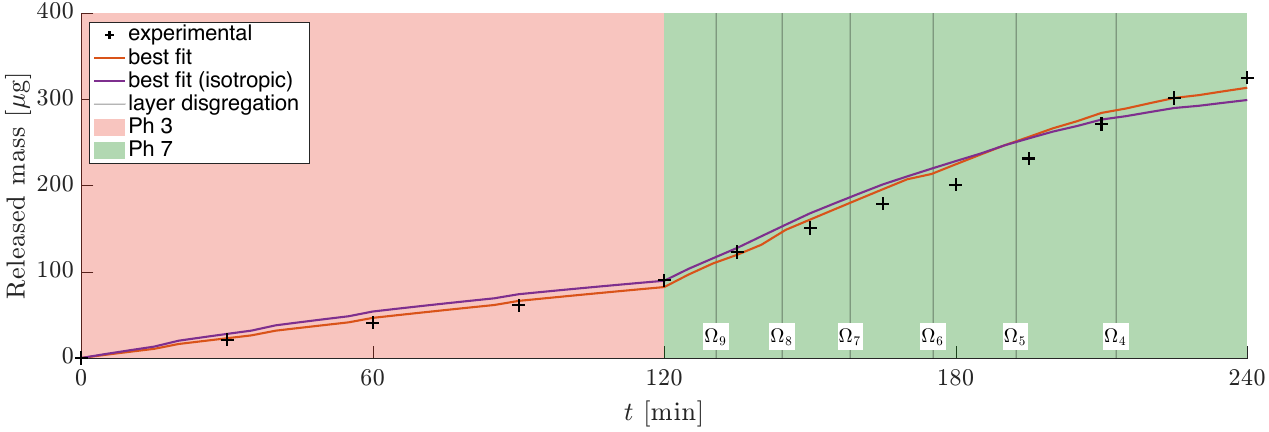}
    \caption{Simulated released TA mass as a function of time, compared with experimental data (+). Vertical lines show the  stratum breakup time points as a consequence of the erosion. The slope of the profile changes after 120 min, related to the faster erosion (cf. fig. 7). A radial anisotropy provides a better fit in the later points.
    }
    \label{fig:experiment-release}
\end{figure}

In \cref{fig:experiment-release}, the experimental data are compared with the predictions of the numerical model, with the TA release curve displayed up to 240 minutes.
After 120 minutes, the rising of pH (due to transition from gastric to intestinal environment) results in a faster release as expected, being yielded by the increase in coating erosion (cf.\ \cref{fig:experiment-profile} and \cref{eq:rs-release-complete}).
Particularly interesting is the shape of the release once reaching $\Omega_3$, which remains nearly constant despite a considerable slowing in the erosion process (cf.\ pre-coting, see also \cref{fig:experiment-profile}).
To carefully reproduce this final experimental change, the core must act as a reservoir of TA, as discussed above, actually requiring the anisotropy of the material for the TA to be released smoothly.
In this sense, for completeness, we include the results of the same computational experiment if an isotropic material is considered instead ($\alpha_\ell = 1, \ell = 1,\dots,L$).
Here, the overall behavioural difference is negligible, but when the pre-coating is reached, where the isotropic experiment actually slows down the TA release due to progressive emptying of the coating itself.

\begin{figure}[t]
    \centering
    \hfill
    \includegraphics[width=.3\linewidth]{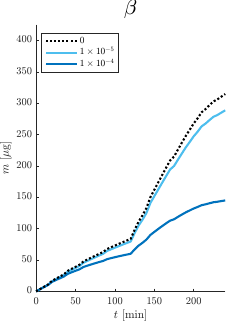}
    \hfill
    \includegraphics[width=.3\linewidth]{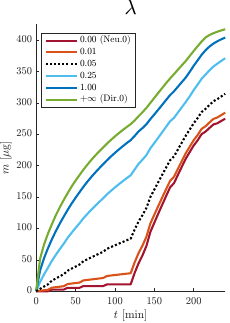}
    \hfill\null
    \\[1em]
    \hfill
    \includegraphics[width=.3\linewidth]{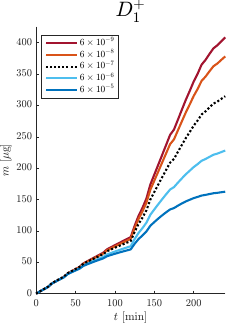}
    \hfill
    \includegraphics[width=.3\linewidth]{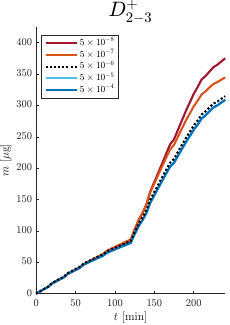}
    \hfill
    \includegraphics[width=.3\linewidth]{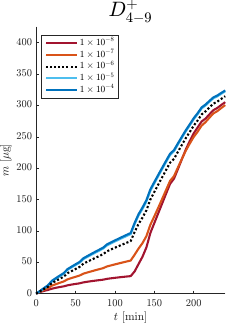}
    \hfill\null
    
    \caption{Sensitivity analysis of parameters $\beta, \lambda, D^+$ on the release curves. The dotted black curve refers to the reference values, the other curves refer to larger/smaller values, as in the legends.}
    \label{fig:sensitivity}
\end{figure}

\paragraph{Parameters sensitivity}

The present model depends on a number of parameters bounded to finite ranges to ensure physical compatibility.
To understand the relevance of such parameters within the simulation environment, we performed a second set of simulations to analyse the sensitivity of $D^+$, $\beta$, $\lambda$, \ie the parameters most contributing to the overall shape of the releasing curve.
In detail, we used the estimation from \cref{ssec:estimating-parameters} as reference values, and we progressively analysed one parameter at a time to explore the model different regimes.
The sensitivity results are presented in \cref{fig:sensitivity}, where the details of the implemented variation (order of magnitude for $D_\ell$, feasible values for $\beta$ and $\lambda$) are detailed within the plots legends.

Unsurprisingly, the effect of the retention (binding or decay) parameter $\beta$ is to reduce the released mass.
A larger mass transfer coefficient $\lambda$ induces a faster release --up to $\lambda = +\infty$ corresponding to a perfect sink-- and smooths out the release curve; conversely, its reduction acts as a barrier to the outgoing mass, where $\lambda = 0$ identifies a perfect barrier allowing no TA diffusion, and consequently isolating the erosion contribution to release.
Core diffusion coefficient $D_1$ positively correlates with TA evasion towards the inner part of the capsule, actually leaving less mass to be released due to erosion, particularly visible in the intestinal phase of the experiment (\ie $t>120$ min).
Analogous is the outcome delivered by the pre-coating diffusion coefficients $D_2$ and $D_3$, whose reduction cause less mass to evade towards the core, actually enhancing the release in the intestinal phase; do notice that a value increase does not correspond to a release increase, suggesting a stable point for the dynamic is reached.
Finally, conversely to the others, outer coating diffusion coefficients ($D_\ell$, $\ell=4, \dots, 9$) positively correlate with mass release, as expected in regular diffusive processes, allowing for higher release of mass; their contribution is however more limited in the intestinal phase, where erosion leads the release process.

\section{Conclusions and future work}

Despite significant advances in the development of innovative and smart drug delivery systems, the design and manufacturing of functionalised carriers --particularly those involving multi-stratum architectures-- remain an open challenge.
The release of therapeutic agents from such systems is governed by a complex interplay of chemical, physical, and structural phenomena that are not yet fully understood.
In this context, mathematical and computational modelling has become a powerful complementary tool, enabling researchers to probe key mechanisms, reduce the need for extensive experimental testing, and explore diverse scenarios \insilico.

In this study, we have introduced a mechanistic model to describe drug release from nano-functio\-na\-lised spherical microcapsules with multi-strata polymer coatings.
The model captures coupled diffusive and reactive processes across a multi-scale structure, where micrometric-scale capsules are surrounded by nanometric strata.
Modelling transport phenomena across such heterogeneous domains presented significant challenges, particularly in handling variable anisotropic properties and interfacial coupling.

Grounded in an experimental setup involving a standardised gastrointestinal simulation protocol (INFOGEST 2.0), the proposed model successfully reproduced the observed \invitro release profiles.
In particular, it captured the releasing effect of the progressive reduction in capsule volume due to erosion, as well as the distinct release phases.
We also found that anisotropy in the early pre-coating strata is essential for correctly describing the dynamics observed during the final release phase, suggesting that microstructural features might affect macroscopic behaviour.
These findings reinforce the relevance of modelling as a tool to interpret complex datasets and guide formulation choices and demonstrate the importance of tightly integrating modelling and experimentation to extract mechanistic insight and validate predictive capabilities.
Concurrently, our computational strategy also facilitates the rational exploration of delivery systems that can meet therapeutic targets with reduced reliance on extensive experimental screening.

In future work, we plan to extend the applicability of the model by including capsule full disintegration and co-delivery of multiple compounds, as well as potentially irregular geometries and heterogeneous coatings with different strata designs.
Additionally, simulating patient-specific or pathological gastrointestinal environments (\eg, inflammation, ulcers, or altered enzyme activity) would improve the model translational relevance.
Coupling the current framework with pharmacokinetic, pharmacodynamics, immune response models, or co-delivery scenarios might further enhance the model utility for pre-clinical designing and optimisation of next-generation drug delivery platforms.

\section*{Authors' contribution}

Conceptualisation: EO, EC, JGH, PG, GP.
Data curation: JGH, AM, PG.
Formal analysis: EO, EC.
Funding acquisition: JGH, PG, GP.
Investigation -- Numerical: EO.
Investigation -- Experimental: JGH, AM, PG.
Methodology -- Numerical: EO, EC.
Methodology -- Experimental: JGH, PG.
Resources: JGH, AM, PG.
Software: EO.
Validation: EO.
Visualization: EO.
Writing -- original draft: EO, EC, JGH, PG, GP.
Writing -- review \& editing: EO, EC, GP.
\\
All authors have read and agreed to the published version of the manuscript.

\section*{Acknowledgements}

This work was supported by 
the NERC Cross--disciplinary Research for Discovery Sciences (NE/ X018229/1), 
the Royal Society International Exchanges 2023 Cost Share Italy (IEC$\backslash$R2$\backslash$232063), 
the Joint Bilateral Agreement CNR (IT) - Royal Society of London (UK) 2024-2025,
the Marie Curie programme (HORIZON-MSCA-2024-PF-01-01 n.\ 101206930-SMartGut Grant), 
and the ATRAE Talento Consolidado programme (ATR2023-145518) funded by the Spanish Ministerio de Ciencia, Innovación y Universidades and the Agencia Estatal de Investigación.

EO and EC are members of the Gruppo Nazionale Calcolo Scientifico -- Istituto Nazionale di Alta Matematica (GNCS--INdAM).

\bibliographystyle{elsarticle-num}
\bibliography{biblio}

\end{document}

%% file: figure/layerByLayer.tex
\begin{tikzpicture}[scale=0.8]
    \def\Radii{{1.8, 2.9, 4.4, 5.0, 5.6, 6.8}}  
    \def\L{5} 
    \def\mycolours{{"blueM","orangeM","yellowM","purpleM","greenM","cyanM"}}
    \def\mylabels{{"1","2","3","","","L"}}

    \foreach \i [
        evaluate=\i as \usecolor using {\mycolours[\i]},
    ] in {5,4,3,2,1,0} {
        \pgfmathsetmacro{\R}{\Radii[\i]}
        \draw[thick, \usecolor, fill=\usecolor!50] (0,0) circle (\R);
    }
    \foreach \i [
        evaluate=\i as \uselabel using {\mylabels[\i]},
    ] in {0,1,2,5} {
        \pgfmathsetmacro{\R}{\Radii[\i]}
        \node[anchor=east] at (\R,0) {\large $\Omega_{\uselabel}$};
        \pgfmathsetmacro{\angle}{(\i+5) * 360/15} 
        \draw[thick, -latex] (0,0) -- ({\R*cos(\angle)}, {\R*sin(\angle)}) node[left]{$R_{\uselabel}$};
    }
    {
        \pgfmathsetmacro{\R}{\Radii[3]}
        \node[] at (\R,0) {\large $\dots$};
    }
    {
        \pgfmathsetmacro{\R}{\Radii[\L]}
        \draw[ultra thick, dashed] (0, 0) circle (\R);
        \node[anchor=south] at (0, \R) {Membrane};
    }
    \draw[latex-latex, thick] (-2.9, 0) --node[below]{$\Delta R_3$} (-4.4, 0.0); 
    \draw[fill=black] (0, 0) circle (0.1);
    \node[] () at (0.7, 0.9) {Core};
\end{tikzpicture}

%% file: figure/layerRadialConstruction.tex
\begin{tikzpicture}[scale=1.35, every node/.style={font=\small}]

    \def\a{0}
    \def\b{1}
    \def\c{2}
    \def\d{3}
    \def\e{6}
    \def\f{7}
    \def\g{8}
    \def\h{9}
    
    \coordinate (A) at (0,0);
    \coordinate (B) at (1,0);
    \coordinate (C) at (2,0);
    \coordinate (D) at (3,0);
    \coordinate (E) at (4,0);
    \coordinate (F) at (5,0);
    \coordinate (G) at (6,0);
    \coordinate (H) at (7,0);
    
    \draw[thick, ->] (\a, 0) -- ($(\h+.5,0)$);
    
    \draw[line width=1pt] (\a, -.1) -- ++(0,0.2);
    \foreach \x in {\b, \d, \e, \g, \h}
        \draw[thick, double distance=.05em] (\x, -.1) -- ++(0,0.2);
    
    \node[below] at (\a,-.1) {$\mathcal{O}\equiv R_0$};
    \node[below] at (\b,-.1) {$R_1$};
    \node[fill=white] at (\c,0)   {$\cdots$};
    \node[below] at (\d,-.1) {$R_{\ell-1}$};
    \node[below] at (\e,-.1) {$R_\ell$};
    \node[fill=white] at (\f,0)   {$\cdots$};
    \node[below] at (\g,-.1) {$R_{L-1}$};
    \node[below] at (\h,-.1) {$R_L$};
    
    \draw[blueM, thick, rounded corners=.6em] 
    (2.2,0.3) -- (2.6,0.4) -- (3.0,0.6) -- (3.4,1.2) -- (3.8,1.0) -- (4.2,0.9) -- (4.6,0.7) -- (5.0,0.4) -- (5.5, 0.4) -- (6.0, 0.35) -- (6.4, 0.2) -- (6.8,0.1);
    \draw[blueM, thick, dashed] (2.98,0.63) -- (2.98,0.1);
    \draw[blueM, thick, dashed] (3.02,0.667) -- (3.02,0.1);
    \draw[white, line width = 0.05em] (3.0,-0.15) -- (3.0,0.7);
    \draw[blueM, thick, dashed] (5.98,0.345) -- (5.98,0.1);
    \draw[blueM, thick, dashed] (6.02,0.335) -- (6.02,0.1);
    \draw[white, line width = 0.05em] (6.0,-0.15) -- (6.0,0.4);

    \node[blueM] at (4.8,0.9) {$c$};
    
    \draw[orangeM] (2.20,0.3) --node[yshift=-.5em] {\tiny$D^-_{\ell-1}$} (2.98,0.3);
    \draw[orangeM] (3.02,0.5) --node[yshift=-.5em] {\tiny$D^-_{\ell}$} (3.48,0.5);
    \draw[orangeM] (3.48,0.7) --node[yshift=-.5em, xshift=-.8em] {\tiny$D^+_{\ell}$} (5.98,0.7);
    \draw[orangeM] (6.02,0.25) --node[yshift=+.5em] {\tiny$D^+_{\ell+1}$} (6.8,0.25);
    \draw[orangeM, dotted] (3.48,0.7) -- (3.48,0.5);
    \node[fill=orangeM, circle, inner sep=1.5pt] at (3.0, 0.4) {};
    \node[fill=orangeM, circle, inner sep=1.5pt] at (6.0, 0.475) {};
    
    
    \node[fill=greenM, circle, inner sep=1.5pt] at (3.48,1.25) {};
    \node[color=greenM, above] at (3.48,1.25) {$F$};
    \draw[<-, thick, greenM] (4.22,1.25) -- (3.48,1.25); 
    \draw[<-, thick, greenM] (3.02,1.25) -- (3.48,1.25);
    
    \fill[red, opacity=0.3] (\h,-.1) rectangle ++(-1.2, .2);
    \draw[red, thick, pattern=north west lines, pattern color=red] (\h,-.1) rectangle ++(-1.2, .2);
    \draw[->, red, thick] (7.8,0.05) --node[above]{$v(t)$} ++(-0.8,0.0);

\end{tikzpicture}

%% file: figure/numericScheme.tex
\begin{tikzpicture}[scale=1.2, every node/.style={font=\small}]

    \def\a{-1}
    \def\b{5} 
    \def\c{9.5}

    \draw[-, blueM, thick] (\a,0) -- (\b,0);
    \foreach \x in {0.25,1.5,2.75,3.9}
        \draw[-, blueM, thick] (\x,-.2) -- (\x,.2);
    \node[blueM] at (0.9, .5) {$\Delta r_{\ell-1}$}; 
    \node at (1.5, -0.7) {$\Omega_{\ell-1}$}; 

    \draw[-, greenM, thick] (\b,0) -- (\c,0);
    \foreach \x in {4.5,5,...,9.5}
        \draw[-, greenM, thick] (\x,-.2) -- (\x,.2);
    \draw[dashed, greenM, thick] (4.5,.15) -- (\b,.15);
    
    \node[greenM] at (6.75, .5) {$\Delta r_{\ell}$};
    \node at (6.75, -0.7) {$\Omega_{\ell}$}; 
    
    \node[greenM] at (6.25, -.3) {\footnotesize$j$-1};
    \node[greenM] at (6.75, -.3) {\footnotesize$j$};
    \node[greenM] at (7.25, -.3) {\footnotesize$j$+1};

    \draw[-, thick] (\a, .5) -- (\a, -.5) node[below] {$R_{\ell-2}$};
    \draw[-, thick] (\b, .5) -- (\b, -.5) node[below] {$R_{\ell-1}$};
    \draw[-, thick] (\c, .5) -- (\c, -.5) node[below] {$R_{\ell}$};

    
    \def\A{3}
    \def\AB{4}
    \def\B{5}
    \def\BC{6}
    \def\C{7}
    \def\CD{8}
    \def\D{9}

    \foreach \x in {\A, \B, \C, \D}
        \draw[-, thick] (\x,2.2) -- (\x,2.8);

    \draw[-, thick] (2.7, 2.5) -- (9.3, 2.5);
    \node[] at (\AB, 2.75) {$C_{j-1}$};
    \node[] at (\AB, 2.495) {$\bullet$};
    \node[] at (\BC, 2.75) {$C_{j}$};
    \node[] at (\BC, 2.495) {$\bullet$};
    \node[] at (\CD, 2.75) {$C_{j+1}$};
    \node[] at (\CD, 2.495) {$\bullet$};
    \node[] at (\AB, 2.2) {$r_{j-1}$};
    \node[] at (\BC, 2.2) {$r_{j}$};
    \node[] at (\CD, 2.2) {$r_{j+1}$};
    \node[] at (5.2, 1.8) {$\hat{D}^\ell_{j-\sfrac12}$};
    \node[] at (7.2, 1.8) {$\hat{D}^\ell_{j+\sfrac12}$};
    \node[] at (5.2, 3.1) {${F}_{j-\sfrac12}$};
    \node[] at (7.2, 3.1) {${F}_{j+\sfrac12}$};

    \draw[-, dashed, redM] (\A,2.15) -- (6, .25);
    \draw[-, dashed, redM] (\D,2.15) -- (7.5, .25);

\end{tikzpicture}